\documentclass[reqno,a4paper,11pt]{article}
\usepackage{amsmath,amsthm,dsfont,amsfonts,amssymb,fancyhdr}
\usepackage{enumerate}
\usepackage[usenames,dvipsnames]{color}
\usepackage{bbm,soul,supertabular,longtable,verbatim,extarrows}
\usepackage{titlesec}
\usepackage{graphicx}
\usepackage[titletoc,toc,title]{appendix}

\usepackage{tikz}
\usetikzlibrary{arrows,decorations.pathmorphing,backgrounds,positioning,fit,petri}

\usepackage{caption}

\usepackage{float}
\usepackage{BOONDOX-cal}
\usepackage{tikz}
\usepackage{booktabs,multirow,makecell}
\usepackage{authblk}
\usepackage[titletoc]{appendix}
\usetikzlibrary[trees]

\usepackage{lipsum}
\usepackage{mathrsfs}
\usepackage{indentfirst}
\usepackage[colorlinks=true, allcolors=blue]{hyperref}
\usepackage[top=2.4cm,bottom=2.2cm,left=2.6cm,right=2cm]{geometry}

\newtheorem{thm}{Theorem}[section]
\newtheorem{lem}[thm]{Lemma}

\newtheorem{de}[thm]{Definition}

\newtheorem{cor}[thm]{Corollary}

\newtheorem{conjecture}[thm]{Conjecture}

\soulregister\cite7
\soulregister\citep7
\soulregister\citet7
\soulregister\ref7
\soulregister\pageref7
\sethlcolor{yellow}
\setstcolor{red}
\soulregister{\em}{0}

\usepackage{xcolor}
\definecolor{backgroundcolor}{RGB}{199, 237, 204}

\begin{document}
	\title{A structure theory for regular graphs with fixed smallest eigenvalue}
	\author[a,b]{Qianqian Yang}
	\author[c,d]{Jack H. Koolen\footnote{J.H. Koolen is the Corresponding author.}}
	\affil[a] {\footnotesize{Department of Mathematics, Shanghai University, Shanghai 200444, PR China}}
	\affil[b] {\footnotesize{Newtouch Center for Mathematics of Shanghai University, Shanghai 200444, PR China}}
	\affil[c]{\footnotesize{School of Mathematical Sciences, University of Science and Technology of China, Hefei, Anhui 230026, PR China}}
	\affil[d]{\footnotesize{CAS Wu Wen-Tsun Key Laboratory of Mathematics, University of Science and Technology of China, Hefei, Anhui 230026, PR China}}
	
	\maketitle
	\pagestyle{plain}
	
	\newcommand\blfootnote[1]{%
		\begingroup
		\renewcommand\thefootnote{}\footnote{#1}%
		\addtocounter{footnote}{-1}%
		\endgroup}
	\blfootnote{E-mail addresses:  {\tt qqyang@shu.edu.cn} (Q. Yang),  {\tt koolen@ustc.edu.cn} (J.H. Koolen).}

\begin{abstract}
In this paper we will give a structure theory for regular graphs with fixed smallest eigenvalue. As a consequence of this theory, we show that a $k$-regular graph with smallest eigenvalue $-\lambda$ has clique number linear in $k$ if $k$ is large with respect to $\lambda$. 
\end{abstract}


\section{Introduction}

In this paper, all graphs are finite, undirected and simple. Suppose that $G$ is a graph. We always denote by $V(G)$ its vertex set and $E(G)$ its edge set. Let $x$ be a vertex of $G$. We also denote by $N_G(x)$ the set of neighbors of $x$ in $G$. For the eigenvalues of $G$, we mean the eigenvalues of its adjacency matrix, and we denote its smallest eigenvalue by $\lambda_{\min}(G)$.

In \cite{ivanov2023} we gave some results without proof. In this paper we supply the proofs for these results. Let $H(a,m)$ be the graph on $a+m+1$ vertices consisting of a clique $K_{a+m}$ together with a vertex which is adjacent to precisely $a$ vertices of this clique. If $a=m$, we also denote $H(m,m)$ by $\widetilde{K}_{2m}$.

In \cite{Greaves.2021}, Greaves, Koolen, and Park showed the following.

\begin{lem}\label{lem:hat}
	Let $G$ be a graph with $\lambda_{\min}(G)\geq-\lambda$ that contains $H(a,m)$ as an induced subgraph.
	Then the inequlity
	\begin{equation}
		\label{eqn:atbound}
		(a-\lambda(\lambda-1))(m-(\lambda-1)^2) \leqslant (\lambda(\lambda-1))^2
	\end{equation}
 holds.
\end{lem}

As a consequence, we have the following results.

\begin{lem}\label{Hat}
	Let $\lambda\geq1$ be an integer and $G$ a graph with $\lambda_{\min}(G) \geq - \lambda$. Let $C$ be a clique in $G$ with order $c$. If $c\geq(\lambda(\lambda-1))^2+\lambda(\lambda-1)+(\lambda-1)^2+2= \lambda^4-2\lambda^3+3\lambda^2-3\lambda+3$, then every vertex outside $C$ has either at most $\lambda(\lambda-1)$ neighbors in $C$, or at least $c-(\lambda-1)^2$ neighbors in $C$.
\end{lem}

Lemma \ref{Hat} can be generalized to $t$-plexes. Let $t$ be a positive integer. An induced subgraph is called a \emph{$t$-plex} if every vertex in it is adjacent to all but at most $t-1$ of the other vertices. Note that a clique is exactly a $1$-plex.

\begin{lem}\label{Hoffproperties1}
	Let $\lambda\geq1$ be an integer and $G$ a graph with  $\lambda_{\min}(G) \geq - \lambda$. Let $P$ be a $t$-plex in $G$ with order $p$. If $p\geq t((\lambda(\lambda-1))^2+\lambda(\lambda-1)+(\lambda-1)^2+(t-1)((\lambda-1)^2+1)+2)$, then every vertex outside $P$ has either at most $t\lambda(\lambda-1)$ neighbors in $P$, or at least $p-t(\lambda-1)^2$ neighbors in $P$.
\end{lem}

In this paper, we will show that a $k$-regular graph with smallest eigenvalue at least $-\lambda$ has a nice family of $((\lambda-1)^2+1)$-plexes if $k$ is large compared to $\lambda$.

\begin{thm}\label{main regular graphs}
	Let $\lambda\geq1$ be an integer. There exists a positive integer $n_1(\lambda)\geq((\lambda-1)^2+1)((\lambda(\lambda-1))^2+\lambda(\lambda-1)+(\lambda-1)^2+(\lambda-1)^2((\lambda-1)^2+1)+2)$ such that for any integer $n\geq n_1(\lambda)$, two integers $K_1(\lambda,n)$ and $C_1(\lambda,n)$ satisfying the following exist.
	
	For any connected $k$-regular graph $G$ with $\lambda_{\min}(G)\geq-\lambda$, if $k\geq K_1(\lambda,n)$, then there exist (maximal) $((\lambda-1)^2+1)$-plexes $P_1, P_2,\ldots, P_{s_1}$ in $G$ satisfying the following five properties:
	\begin{enumerate}
		\item  $1\leq\{i=1,2,...,s_1 \mid x\in V(P_i)\}|\leq \lambda $ for every vertex $x$ of $G$.
		\item  For every vertex $x$ of $G$, $|\{y\mid y\text{ is adjacent to }x\}-\{y\mid x,y\in V(P_i)\text{ for some }i=1,...,s_1\}|\leq C_1(\lambda,n)$ holds.
		\item $|V(P_i)|\geq n$ for $i=1,\ldots,s_1$.
		\item  $|V(P_i)\cap V(P_j)|\leq\lambda-1$ for any $1\leq i<j\leq s_1$.
		\item Every vertex not in $P_i$ $(i=1,2,\ldots,s_1)$ has at most $((\lambda-1)^2+1)\lambda(\lambda-1)$ neighbors in $P_i$.
	\end{enumerate}
\end{thm}

As an application of Theorem \ref{main regular graphs}, we obtain that

\begin{cor}\label{clique size}
	Let $\lambda\geq1$ be an integer. There exists a positive integer $\kappa(\lambda)$ and a positive real number $\varepsilon(\lambda)$ such that for any connected $k$-regular graph $G$ with $\lambda_{\min}(G)\geq-\lambda$, if $k\geq \kappa(\lambda)$, then every vertex lies in a clique of order at least $\frac{k-\varepsilon(\lambda)}{\lambda((\lambda-1)^2+1)}$.
\end{cor}

Let $\lambda\geq2$ be an integer. Given a graph $G$, let $\omega(G)$ be the clique number of the graph $G$.
Define 
\[
\gamma(\lambda):=\lim_{k\to\infty}\inf\left\{\frac{\omega(G)}{k} \mid G \text{ is a connected $k$-regular graph and }\lambda_{\min}(G) \geq -\lambda\right\}.
\]

The Hamming graphs $H(D, q)$ have valency
$D(q-1)$, smallest eigenvalue $-D$,  and $\omega(H(D, q))=q$. This shows that 
\[
\gamma(D) \leq \frac{1}{D} \text{ for any integer }D\geq2.
\]
Note that 
\[
\gamma(2) =\frac{1}{2},
\] 
as for a connected regular graph $G$ with at least $29$ vertices, either $G$ is a line graph of a regular or bipartite and semiregular graph, or $G$ is a cocktail party graph (see \cite[Theorem 4.4, Theorem 4.11]{Cameron} and \cite[Lemma 1.7.5]{GD01}). 

In 1995, Woo and Neumaier \cite[p.~590, Conjecture (\rm{ii})]{Woo} conjectured that a regular graph with smallest eigenvalue at least $-1-\sqrt{2}$ and sufficiently large valency is a line graph or a cocktail party graph. However Yu \cite{Yu}, by using Hoffman graphs, constructed an infinite family of regular graphs with smallest eigenvalue in $[-1-\sqrt{2},-2)$ and disproved this conjecture later. We use now a similar method to construct a family of connected regular graphs with smallest eigenvalue at least $-3$ in which every vertex lies in a clique of order approximately a quarter of the valency and show that
\begin{equation}\label{gamm3 upper bound}
	\gamma(3) \leq \frac{1}{4}.
\end{equation}
The construction relies on the Steiner systems $S(2,4,v)$.

A Steiner system $S(2,4,v)$ is a pair $(X,\mathfrak{B})$, where $X$ is a $v$-element set and $\mathfrak{B}$ is a family of $4$-element subsets of $X$ called \emph{blocks} such that each $2$-element subset of $X$ is contained in exactly one block. A Steiner system $S(2,4,v)$ exists if and only if $v\equiv 1$ or $4\pmod {12}$ (see \cite{Hanani1961}). Let $(X,\mathfrak{B})$ be a Steiner system and let $\mathfrak{s}=\uplus_{i=1}^b\mathfrak{s}_i$ be a Hoffman graph, where $\mathfrak{s}_i$ is isomorphic to the following Hoffman graph

	\begin{center}
	\raisebox{0.1ex}
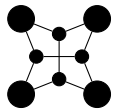
{\begin{tikzpicture}[scale=1]
			\draw (0,-0.8) node {$$};
			\tikzstyle{every node}=[draw,circle,fill=black,minimum size=10pt,scale=1,
			inner sep=0pt]

			\draw (-0.5,0) node (2f1) [label=below:$$] {};
			\draw (0.5,0) node (2f2) [label=below:$$] {};
			\draw (-0.5,1) node (2f3) [label=below:$$] {};
			\draw (0.5,1) node (2f4) [label=below:$$] {};

			\tikzstyle{every node}=[draw,circle,fill=black,minimum size=5pt,scale=1,
			inner sep=0pt]

			\draw (0,0.2) node (2s1) [label=below:$$] {};
			\draw (0.3,0.5) node (2s2) [label=below:$$] {};
			\draw (-0.3,0.5) node (2s3) [label=below:$$] {};
			\draw (0,0.8) node (2s4) [label=below:$$] {};

			\draw (2f1) -- (2s1) -- (2f2) -- (2s2) -- (2f4) -- (2s4) -- (2f3) -- (2s3) -- (2f1);
			\draw (2s1) -- (2s4);
			\draw (2s2) -- (2s3);
	\end{tikzpicture}}
\end{center}
\vspace{-1cm}
for $i=1,\ldots,b$, $V_{\mathrm{fat}}(\mathfrak{s})=X$ and $\{V_{\mathrm{fat}}(\mathfrak{s}_1),\ldots,V_{\mathrm{fat}}(\mathfrak{s}_b)\}=\mathfrak{B}$. Note that for each fat vertex of $\mathfrak{s}$, there are eaxactly $r=\frac{v-1}{3}$ Hoffman graphs among $\mathfrak{s}_i,\ldots,\mathfrak{s}_b$ having it as fat vertex. So, the slim graph of $\mathfrak{s}$ is a connected $(4r-3)$-regular graph and every vertex in it lies in a clique of order $r$.  This shows the inequality \eqref{gamm3 upper bound}.

Corollary \ref{clique size} implies that 
\[
\gamma(\lambda) \geq \frac{1}{\lambda(\lambda-1)^2+1}.
\]

We think our bound on  $\gamma(\lambda)$ can be improved and we conjecture that

\begin{conjecture}
There exist positive constants $\varepsilon_1$ and $\varepsilon_2$ such that 
$\frac{\varepsilon_1}{\lambda^2}\leq \gamma(\lambda)\leq\frac{\varepsilon_2}{\lambda^2}$ holds for every integer $\lambda$.
\end{conjecture}

This paper is organized as follows. In Section \ref{sec:Hoffman} we will introduce our main tool: Hoffman graphs, and in Section \ref{sec:proof}, we will show the proofs of our results mentioned in Introduction.

\section{Preliminaries}\label{sec:Hoffman}
\subsection{Hoffman graphs}
\begin{de} A \emph{Hoffman graph} $\mathfrak{h}$ is a pair $(H,\ell)$, where $H=(V,E)$ is a graph and $\ell:V \to\{f,s\}$ is a labeling map satisfying the following conditions:
	\begin{enumerate}
		\item vertices with label $f$ are pairwise non-adjacent,
		\item every vertex with label $f$ is adjacent to at least one vertex with label $s$.
	\end{enumerate}
\end{de}
We call a vertex with label $s$ a \emph{slim vertex}, and a vertex with label $f$ a \emph{fat vertex}. We denote by $V_{\mathrm{slim}}(\mathfrak{h})$ (resp. $V_{\mathrm{fat}}(\mathfrak{h})$) the set of slim (resp. fat) vertices of $\mathfrak{h}$.

\vspace{0.1cm}
For a vertex $x$ of $\mathfrak{h}$, we denote by $N_{\mathfrak{h}}^{s}(x)$ (resp. $N_{\mathfrak{h}}^{f}(x)$) the set of slim (resp. fat) neighbors of $x$ in $\mathfrak{h}$. If every slim vertex of $\mathfrak{h}$ has a fat neighbor, then we call $\mathfrak{h}$ \emph{fat}, and if every slim vertex of $\mathfrak{h}$ has at least $t$ fat neighbors, we call $\mathfrak{h}$ $t$-\emph{fat}. In a similar fashion, we denote by $N^{f}_\mathfrak{h}(x_1,x_2)$ the set of common fat neighbors of two slim vertices $x_1$ and $x_2$ in $\mathfrak{h}$ and $N^{s}_\mathfrak{h}(f_1,f_2)$ the set of common slim neighbors of two fat vertices $f_1$ and $f_2$ in $\mathfrak{h}$.

\vspace{0.1cm}
The \emph{slim graph} of the Hoffman graph $\mathfrak{h}$ is the subgraph of $H$ induced on $V_{\mathrm{slim}}(\mathfrak{h})$. Note that any graph can be considered as a Hoffman graph with only slim vertices, and vice versa. We will not distinguish between Hoffman graphs with only slim vertices and graphs.

\vspace{0.1cm}
A Hoffman graph $\mathfrak{h}_1= (H_1, \ell_1)$ is called an (\emph{proper}) \emph{induced Hoffman subgraph} of $\mathfrak{h}=(H, \ell)$, if $H_1$ is an (proper) induced subgraph of $H$ and $\ell_1(x) = \ell(x)$ holds for all vertices $x$ of $H_1$.

\vspace{0.1cm}
Let $W$ be a subset of $V_{\mathrm{slim}}(\mathfrak{h})$. An induced Hoffman subgraph of $\mathfrak{h}$ generated by $W$, denoted by $\langle W\rangle_{\mathfrak{h}}$, is the Hoffman subgraph of $\mathfrak{h}$ induced on $W \cup\{f\in V_{\mathrm{fat}}(\mathfrak{h})\mid f \sim w \text{ for some }w\in W \}$.

\vspace{0.1cm}
For a fat vertex $f$ of $\mathfrak{h}$, a \emph{quasi-clique} (with respect to $f$) is a subgraph of the slim graph of $\mathfrak{h}$ induced on the slim vertices adjacent to $f$ in $\mathfrak{h}$, and we denote it by $Q_{\mathfrak{h}}(f)$.
\begin{de}Two Hoffman graphs $\mathfrak{h}= (H, \ell)$ and $\mathfrak{h}^\prime=(H^\prime, \ell^\prime)$ are \emph{isomorphic} if there exists an isomorphism from $H$ to $H^\prime$ which preserves the labeling.
\end{de}



\begin{de} For a Hoffman graph $\mathfrak{h}=(H,\ell)$, there exists a matrix $C$ such that the adjacency matrix $A$ of $H$ satisfies
	\begin{eqnarray*}
		A=\left(
		\begin{array}{cc}
			A_s  & C\\
			C^{T}  & O
		\end{array}
		\right),
	\end{eqnarray*}
	where $A_s$ is the adjacency matrix of the slim graph of $\mathfrak{h}$. The \emph{special matrix} $Sp(\mathfrak{h})$ of $\mathfrak{h}$ is the real symmetric matrix $A_s-CC^{T}.$
\end{de}

The \emph{eigenvalues} of $\mathfrak{h}$ are the eigenvalues of its special matrix $Sp(\mathfrak{h})$, and the smallest eigenvalue of $\mathfrak{h}$ is denoted by  $\lambda_{\min}(\mathfrak{h})$. Note that $\mathfrak{h}$ is not determined by its special matrix in general, since different $\mathfrak{h}$'s may have the same special matrix. Observe also that if there are no fat vertices in $\mathfrak{h}$, then $Sp(\mathfrak{h})=A_s$ is just the standard adjacency matrix.

\begin{lem}[{\cite[Lemma 3.4]{Woo}}]\label{fatnbr}
	Let $\mathfrak{h}$ be a Hoffman graph and let $x_i$ and $x_j$ be two distinct slim vertices of $\mathfrak{h}$. The special matrix $Sp(\mathfrak{h})$ has diagonal entries
	
	$$Sp(\mathfrak{h})_{x_i,x_i} = - |N_\mathfrak{h}^f (x_i)|$$
	and off-diagonal entries
	$$Sp(\mathfrak{h})_{x_i,x_j} =\left\{\begin{array}{ll}
	1-|N_\mathfrak{h}^f (x_i,x_j)|,& x_i,x_j\text{~are adjacent},\\
	-|N_\mathfrak{h}^f (x_i,x_j)|, &x_i,x_j\text{~are not adjacent}.
	\end{array}
	\right. $$
\end{lem}

For the smallest eigenvalues of Hoffman graphs and their induced Hoffman subgraphs, Woo and Neumaier showed the following inequality.
\begin{lem}[{\cite[Corollary 3.3]{Woo}}]\label{hoff}
	If $\mathfrak{h}_1$ is an induced Hoffman subgraph of a Hoffman graph $\mathfrak{h}$, then $\lambda_{\min}(\mathfrak{h}_1)\geq\lambda_{\min}(\mathfrak{h})$ holds.
\end{lem}

As a corollary of Lemma \ref{hoff}, we have
\begin{lem}\label{interelacing}
	If $G_1$ is an induced subgraph of $G$, then $\lambda_{\min}(G_1)\geq\lambda_{\min}(G)$ holds.
\end{lem}


%

\begin{de}\label{directsummatrix}
	Let $\mathfrak{h}^1$ and $\mathfrak{h}^2$ be two Hoffman graphs. A Hoffman graph $\mathfrak{h}$ is the \emph{sum} of $\mathfrak{h}^1$ and $\mathfrak{h}^2$, denoted by $\mathfrak{h} =\mathfrak{h}^1\uplus\mathfrak{h}^2$, if $\mathfrak{h}$ satisfies the following condition:
	
	There exists a partition $\big\{V_{\mathrm{slim}}^1(\mathfrak{h}),V_{\mathrm{slim}}^2(\mathfrak{h})\big\}$ of $V_{\mathrm{slim}}(\mathfrak{h})$ such that induced Hoffman subgraphs generated by $V_{\mathrm{slim}}^i(\mathfrak{h})$ are $\mathfrak{h}^i$ for $i=1,2$ and
	\[Sp(\mathfrak{h})=
	\begin{pmatrix}
		Sp(\mathfrak{h}^1) & O \\
		O& Sp(\mathfrak{h}^2)
	\end{pmatrix}
	\] with respect to the partition $\big\{V_{\mathrm{slim}}^1(\mathfrak{h}),V_{\mathrm{slim}}^2(\mathfrak{h})\big\}$ of $V_{\mathrm{slim}}(\mathfrak{h})$.
\end{de}

Clearly, by definition, the sum is associative, so that the sum $\biguplus_{i=1}^{r}\mathfrak{h}^i$ is well-defined. We can check that $\mathfrak{h}$ is a sum of two non-empty Hoffman graphs if and only if $Sp(\mathfrak{h})$ is a block matrix with at least two blocks. If $\mathfrak{h} =\mathfrak{h}^1 \uplus \mathfrak{h}^2$ for some non-empty Hoffman subgraphs $\mathfrak{h}^1$ and $\mathfrak{h}^2$, then we call $\mathfrak{h}$ \emph{decomposable} with $\{\mathfrak{h}^1,\mathfrak{h}^2\}$ as a \emph{decomposition} and call $\mathfrak{h}^1, \mathfrak{h}^2$ \emph{factors} of $\mathfrak{h}$. 
Otherwise, $\mathfrak{h}$ is called \emph{indecomposable}.

The following lemma gives a combinatorial way to define the sum of Hoffman graphs.

\begin{lem}[{\cite[Lemma 2.11]{kyy1}}]\label{combi}
	Let $\mathfrak{h}$ be a Hoffman graph and $\mathfrak{h}^1$ and $\mathfrak{h}^2$ be two induced Hoffman subgraphs of $\mathfrak{h}$. The Hoffman graph $\mathfrak{h}$ is the sum of $\mathfrak{h}^1$ and $\mathfrak{h}^2$ if and only if $\mathfrak{h}^1$, $\mathfrak{h}^2$, and $\mathfrak{h}$ satisfy the following conditions:
	\begin{enumerate}
		\item $V(\mathfrak{h})=V(\mathfrak{h}^1)\cup V(\mathfrak{h}^2);$
		\item $\big\{V_{\mathrm{slim}}(\mathfrak{h}^1),V_{\mathrm{slim}}(\mathfrak{h}^2)\big\}$ is a partition of $V_{\mathrm{slim}}(\mathfrak{h});$
		\item if $x \in V_{\mathrm{slim}}(\mathfrak{h}^i),~f \in V_{\mathrm{fat}}(\mathfrak{h})$ and $x$ is adjacent to $f$ in $\mathfrak{h}$, then $f\in V_{\mathrm{fat}}(\mathfrak{h}^i);$
		\item if $x \in V_{\mathrm{slim}}(\mathfrak{h}^1)$ and $y \in V_{\mathrm{slim}}(\mathfrak{h}^2)$, then $x$ and $y$ have at most one common fat neighbor, and they have one if and only if they are adjacent.
	\end{enumerate}
\end{lem}

\subsection{Graphs and Hoffman graphs}
Let $m$ be a positive integer and $G$ a graph that does not contain $\widetilde{K}_{2m}$ as an induced subgraph, where $\widetilde{K}_{2m}$ is the graph with $2m+1$ vertices consisting of a clique $K_{2m}$ together with a vertex that is adjacent to exactly $m$ vertices of the clique. Let $\mathcal{C}(n)=\{C\mid$ $C$ is a maximal clique of $G$ of order at least $n\}$. Define the relation $\equiv_n^m$ on $\mathcal{C}(n)$ by $C_1\equiv_n^m C_2$ if each vertex $x$ in $C_1$ has at most $m-1$ non-neighbors in $C_2$ and each vertex $y$ in $C_2$ has at most $m-1$ non-neighbors in $C_1$. Note that $\equiv_n^m$ is an equivalence relation if $n \geq (m+1)^2$ (see \cite[Lemma 3.1]{kky}).

Let $[C]_n^m$ denote the equivalence class of $\mathcal{C}(n)$ of $G$ under the equivalence relation $\equiv_n^m$ containing the maximal clique $C$ of $\mathcal{C}(n)$. We define the \emph{quasi-clique} $Q([C]_n^m)$ of $C$ with respect to the pair $(m, n)$ as the subgraph of $G$ induced on the set $\{x\in V(G)\mid$ $x$ has at most $m-1$ non-neighbors in $C\}$. Note that for any $C'\in [C]_n^m$, we have $Q([C']_n^m)=Q([C]_n^m)$ (see \cite[Lemma 3.3]{kky}).

Let $[C_1]_n^m, \dots, [C_r]_n^m$ be the equivalence classes of maximal cliques under $\equiv_n^m$. The \emph{associated Hoffman graph} $\mathfrak{g}=\mathfrak{g}(G,m,n)$ is the Hoffman graph satisfying the following conditions:
\begin{enumerate}
	\item $V_{\mathrm{slim}}(\mathfrak{g}) = V(G)$, $V_{\mathrm{fat}}(\mathfrak{g})=\{f_1, f_2, \dots, f_r\}$;
	\item the slim graph of $\mathfrak{g}$ equals $G$;
	\item for each $i$, the fat vertex $f_i$ is adjacent to exactly all the vertices of $Q([C_i]_n^m)$ for $i=1,2, \dots, r.$
\end{enumerate}

From the above definition of associated Hoffman graphs, we find that for each $i=1,\dots,r$, the quasi-clique $Q([C_i]_n^m)$ of $C_i$ with respect to the pair $(m, n)$ is exactly the quasi-clique $Q_{\mathfrak{g}}(f_i)$ in $\mathfrak{g}$ with respect to the fat vertex $f_i$.

%

\section{Proofs}\label{sec:proof}
In this section, we will prove Lemma \ref{Hat}, Lemma \ref{Hoffproperties1}, Theorem \ref{main regular graphs} and Corollary \ref{clique size}.

\begin{proof}[Proof of Lemma \ref{Hat}]
   Assume that a vertex $x$ outside $C$ has exactly $s$ neighbors in $C$, where $\lambda(\lambda-1)+1\leq s\leq c-(\lambda-1)^2-1$. Then $G$ contains $H(s,c-s)$ as an induced subgraph and by Lemma \ref{lem:hat}, we have $(\lambda(\lambda-1))^2\geq (s-\lambda(\lambda-1))(c-s-(\lambda-1)^2)$. Let $f(s)=(s-\lambda(\lambda-1))(c-s-(\lambda-1)^2)$. Then $(\lambda(\lambda-1))^2\geq f(s)\geq\min\{f(\lambda(\lambda-1)+1),f(c-(\lambda-1)^2)-1)\}=c-\lambda(\lambda-1)-(\lambda-1)^2-1$, that is, $c\leq (\lambda(\lambda-1))^2+\lambda(\lambda-1)+(\lambda-1)^2+1$. This gives a contradiction. Hence, the lemma holds.
\end{proof}

Before proving Lemma \ref{Hoffproperties1}, we need a result of Hajnal and Szemer\'{e}di.

\begin{lem}[cf.~{\cite{Hajnal.1970}}]\label{clique decomposition}
Let $t\geq1$ be an integer and $P$ a $t$-plex with vertex set $V(P)$. There exists a partition $V_1\cup V_2\cdots\cup V_t$ of $V(P)$ such that the induced subgraph on $V_i$ is a clique for $i=1,...,t$ and $||V_i|-|V_j||\leq1$ for $1\leq i,j\leq t$.

\end{lem} 

\begin{proof}[Proof of Lemma \ref{Hoffproperties1}]
Let $V_1\cup V_2\cdots\cup V_t$ be a parition of the vertex set $V(P)$ of $P$, where $V_i$ induces a clique $C_i$ for $i=1,\ldots,t$ and $||V_i|-|V_j||\leq1$ for $i,j=1,...,t$. As $|V(P)|=p\geq t((\lambda(\lambda-1))^2+\lambda(\lambda-1)+(\lambda-1)^2+(t-1)((\lambda-1)^2+1)+2)$, then for each $i$, $C_i$ has order at least $(\lambda(\lambda-1))^2+\lambda(\lambda-1)+(\lambda-1)^2+(t-1)((\lambda-1)^2+1)+2$. For a vertex $x$ outside $P$, By Lemma \ref{Hat}  we have either $|N_G(x)\cap V_i|\leq \lambda(\lambda-1)$ or $|N_G(x)\cap V_i|\geq |V_i|-(\lambda-1)^2$, where $N_G(x)$ is the set of neighbors of $x$ in $G$.

Assume $|N_G(x)\cap V_1|\leq \lambda(\lambda-1)$ and $|N_G(x)\cap V_2|\geq |V_2|-(\lambda-1)^2$. Let $y_1,\ldots,y_{(\lambda-1)^2+1}$ be $(\lambda-1)^2+1$ non-neighbors of $x$ in $V_1$ and let $V_2'=V_2\cap N_G(x)\cap (\bigcap_{i=1}^{(\lambda-1)^2+1}N_G(y_i))$. It is not hard to find $|V_2-N_G(y_{\ell})|\leq t-1$ for $\ell=1,\ldots,(\lambda-1)^2+1$, as $P$ is a $t$-plex. So we have $|V_2'|\geq |V_2|-(\lambda-1)^2-(t-1)((\lambda-1)^2+1)\geq  (\lambda(\lambda-1))^2+\lambda(\lambda-1)+2$. Now $\{y_1,\ldots,y_{(\lambda-1)^2+1}\}\cup V_2'$ induces a clique $C$ of order at least $(\lambda(\lambda-1))^2+\lambda(\lambda-1)+(\lambda-1)^2+2$. But for the vertex $x$, $\lambda(\lambda-1)+1\leq |N_G(x)\cap V(C)|=|V_2'|=|V(C)|-((\lambda-1)^2+1)$. This contradicts Lemma \ref{Hat}. Hence
either $|N_G(x)\cap V_i|\leq \lambda(\lambda-1)$ hold for all $i$, or  $|N_G(x)\cap V_i|\geq |V_i|-(\lambda-1)^2$ hold for all $i$. The desired result follows.
\end{proof}

In preparation for proving Theorem \ref{main regular graphs}, we need some results on associated Hoffman graphs of graphs with fixed smallest eigenvalue.

\begin{thm}[{\cite[p.~12]{kky}}]\label{associated Hoff 1}
	Let $\lambda\geq1$ be a positive integer. Let $\mathfrak{h}^{(\lambda)}$ be the Hoffman graph with exactly $\lambda$ fat vertices and  one slim vertex which is adjacent to all these $\lambda$ fat vertices. Let $m(\lambda)=\min\{m\mid \lambda_{\min}(\widetilde{K}_{2m}) <-\lambda\}$. There exists a postive integer $N_1(\lambda)\geq (m(\lambda)+1)^2$ such that if $G$ is a graph with $\lambda_{\min}(G)\geq-\lambda$, then for any integer $n\geq N_1(\lambda)$, the associated Hoffman graph $\mathfrak{g}(G,m(\lambda),n)$ does not contain $\mathfrak{h}^{(\lambda+1)}$ as an induced Hoffman subgraph.
\end{thm}

\begin{thm}[{\cite[Theorem 4.2]{kyy1}}]\label{associated Hoff 2}
	Let $\lambda\geq1$ be a positive integer. Let $m(\lambda)=\min\{m\mid \lambda_{\min}(\widetilde{K}_{2m}) <-\lambda\}$. There exists a postive integer $N_2(\lambda)\geq (m(\lambda)+1)^2$ such that if $G$ is a graph with $\lambda_{\min}(G)\geq-\lambda$, then for any integer $n\geq N_2(\lambda)$, the quasi-cliques $Q_1,Q_2,\ldots,Q_r$ of $G$ with respect to the pair $(m(\lambda),n)$ satisfy the following conditions:
	\begin{enumerate}
		\item the complement of $Q_i$ has valency at most $(\lambda-1)^2$ for $i=1,\ldots,r$;
		\item the intersection $V(Q_i)\cap V(Q_j)$ contains at most $\lambda-1$ vertices for $1\leq i<j\leq r$.
	\end{enumerate} 
\end{thm}

\begin{thm}\label{main Hoffman regular graph}
	Let $\lambda\geq1$ be a positive integer. Let $m(\lambda)=\min\{m\mid \lambda_{\min}(\widetilde{K}_{2m}) <-\lambda\}$. There exists a postive integer $n_2(\lambda)\geq\max \{(m(\lambda)+1)^2,((\lambda-1)^2+1)((\lambda(\lambda-1))^2+\lambda(\lambda-1)+(\lambda-1)^2+(\lambda-1)^2((\lambda-1)^2+1)+2)\}$ such that for any integer $n\geq n_2(\lambda)$, two real numbers $K_2(\lambda,n)$ and $C_2(\lambda,n)$ satisfying the following exist.
	
	For any connected $k$-regular graph $G$ with $\lambda_{\min}(G)\geq-\lambda$, if $k\geq K_2(\lambda,n)$, then the associated Hoffman graph $\mathfrak{g}:=\mathfrak{g}(G,m(\lambda),n)$ has the following property. 
	\begin{enumerate}
		\item For every $x\in V_{\mathrm{slim}}(\mathfrak{g})$, $1\leq |N_{\mathfrak{g}}^f(x)|\leq\lambda$.
		\item For every $x\in V_{\mathrm{slim}}(\mathfrak{g})$, $|N_{\mathfrak{g}}^s(x)-\bigcup_{f\in N_{\mathfrak{g}}^f(x)}N_{\mathfrak{g}}^s(f)|\leq C_2(\lambda,n)$.
		\item  Every quasi-clique in $\mathfrak{g}$ is a $((\lambda-1)^2+1)$-plex of order at least $n$.
		\item For any two quasi-cliques $Q$ and $Q'$ in $\mathfrak{g}$, $|V(Q)\cap V(Q')|\leq\lambda-1$. 
	
		\item For every $x\in V_{\mathrm{slim}}(\mathfrak{g})$ and every quasi-clique $Q$ in $\mathfrak{g}$, if $x\not\in V(Q)$, then $x$ has at most $((\lambda-1)^2+1)\lambda(\lambda-1)$ neighbors in $Q$.
	\end{enumerate}
\end{thm}
\begin{proof}
Let $n_2(\lambda)=\max\{N_1(\lambda),N_2(\lambda),((\lambda-1)^2+1)((\lambda(\lambda-1))^2+\lambda(\lambda-1)+(\lambda-1)^2+(\lambda-1)^2((\lambda-1)^2+1)+2)\}$, where $N_1(\lambda)$ and $N_2(\lambda)$ are such that Theorem \ref{associated Hoff 1} and Theorem \ref{associated Hoff 2} hold respectively. For every integer $n\geq n_2(\lambda)$, let $K_2(\lambda,n)=R(n-1,\lambda^2+1)$ and $C_2(\lambda,n)=R(n-1,\lambda^2+1)-1$, where $R(~,~)$ denote the Ramsey number of two integers. We will show that the associated Hoffman graph $\mathfrak{g}:=\mathfrak{g}(G,m(\lambda),n)$ satisfies {\rm(i)}--{\rm (v)}.

{\rm (i)} As $\lambda_{\min}(G)\geq-\lambda$, $G$ does not contain the $(\lambda^2+1)$-claw graph $K_{1,\lambda^2+1}$, which has smallest eigenvalue $-\sqrt{\lambda^2+1}$, as an induced subgraph by Lemma \ref{interelacing}. Considering $k\geq K_2(\lambda,n)=R(n-1,\lambda^2+1)$, we find that the neighbors of every vertex of $G$ induces a clique of order at least $n-1$. Thus, every vertex of $G$ lies in a clique of order at least $n$. This implies $|N_{\mathfrak{g}}^f(x)|\geq1$.  By Theorem \ref{associated Hoff 1}, we also have $ |N_{\mathfrak{g}}^f(x)|\leq\lambda$.

{\rm (ii)} Let $T=N_{\mathfrak{g}}^s(x)-\bigcup_{f\in N_{\mathfrak{g}}^f(x)}N_{\mathfrak{g}}^s(f)$. 
It is sufficient to show $|T|\leq C_2(\lambda,n)$. Otherwise, supoose $|T|\geq C_2(\lambda,n)+1=R(n-1,\lambda^2+1)$, and $T$ will induces a clique of order at least $n-1$. By the definition of quasi-cliques, there will be a quasi-clique with respect to the pair $(m(\lambda),n)$ which contains $x$ and a clique of order at least $n-1$ in $T$. This contradicts the definition of $T$. Thus, $|T|\leq C_2(\lambda,n)$.

{\rm (iii)}--{\rm (iv)} They follow from Theorem \ref{associated Hoff 2} directly.

{\rm (v)} To prove it, it is sufficient to show, by {\rm (iii)} and Lemma \ref{Hoffproperties1}, that the case where $x$ has at least $|V(Q)|-((\lambda-1)^2+1)(\lambda-1)^2$ neighbors in $Q$ is not impossible. Suppose not. Let $C$ be a maximal clique of order at least $n$ in $Q$. As $|V(Q)\cap N_G(x)|\geq |V(Q)|-((\lambda-1)^2+1)(\lambda-1)^2$, we have $|V(C)-N_G(x)|\leq |V(Q)-N_G(x)|\leq ((\lambda-1)^2+1)(\lambda-1)^2$. Considering that $C$ has order at least $n$, we find that $x$ has at most $(\lambda-1)^2$ non-neighbors in $C$ by Lemma \ref{Hat}. 

Assume $(\lambda-1)^2\leq m(\lambda)-1$. By the definition of quasi-cliques with respect to the pair $(m(\lambda),n)$, we conclude $x\in V(Q)$, which gives us contradiction. So now we are going to prove that the assumption $(\lambda-1)^2\leq m(\lambda)-1$ holds. To do so, we will show that the graph $\widetilde{K}_{2(\lambda-1)^2}$ has smallest eigenvalue at least $-\lambda$. Note that the graph $\widetilde{K}_{2(\lambda-1)^2}$ has the same smallest eigenvalue as the matrix 
\[
M=\begin{pmatrix}
	0 & (\lambda-1)^2 & 0\\
	1 & (\lambda-1)^2-1 & (\lambda-1)^2\\
	0 & (\lambda-1)^2  & (\lambda-1)^2-1
\end{pmatrix}
\]
does. Since $\det(M+\lambda I)=\lambda^2(\lambda-1)^2\geq 0$, $M$ has smallest eigenvalue at least $-\lambda$, so does $\widetilde{K}_{2(\lambda-1)^2}$. This completes the proof.
\end{proof}

\begin{proof}[Proof of Theorem \ref{main regular graphs}]
	This follows from Theorem \ref{main Hoffman regular graph} directly.
\end{proof}

\begin{proof}[Proof of Corollary \ref{clique size}]
Let $\kappa(\lambda)=K_1(\lambda,n_1(\lambda))$ and $\varepsilon(\lambda)=C_1(\lambda,n_1(\lambda))$, where $n_1(\lambda)$, $K_1(\lambda,~)$ and $C_1(\lambda,~)$ are such that Theorem \ref{main regular graphs} holds. By Theorem \ref{main regular graphs} {\rm (i)} and  {\rm (ii)}, every vertex lies in a $((\lambda-1)^2+1)$-plex of order at least $\frac{k-\varepsilon(\lambda)}{\lambda}$. Applying Lemma \ref{clique decomposition}, the desired result follows.
\end{proof}

\section*{Acknowledgements}
Q. Yang is supported by the Shanghai Sailing Program (No. 23YF1412500).

J. H. Koolen is partially supported by the National Key R. and D. Program of China (No. 2020YFA0713100), the National Natural Science Foundation of China (No. 12071454 and No. 12371339), and the 
Anhui Initiative in Quantum Information Technologies (No. AHY150000).

\bibliographystyle{plain}
\bibliography{YK}

\end{document}